\documentclass{article}
\usepackage{anysize}
\marginsize{2.5cm}{2.5cm}{3cm}{3cm}
\usepackage{authblk}
\usepackage{amsmath}
\usepackage{latexsym}
\usepackage{amsfonts}
\usepackage{theorem}
\usepackage{amssymb}
\usepackage{graphicx}
\DeclareGraphicsExtensions{.eps,.png,.pdf,.jpg}
\usepackage{subfig}
\usepackage{multirow}
\usepackage{color}
\usepackage{fancyvrb}
\usepackage[utf8]{inputenc}
\definecolor{mblue}{rgb}{0,0,1}
\definecolor{mgreen}{rgb}{0.13333,0.5451,0.13333}
\definecolor{mred}{rgb}{0.62745,0.12549,0.94118}
\definecolor{mgrey}{rgb}{0.5,0.5,0.5}
\definecolor{mdarkgrey}{rgb}{0.25,0.25,0.25}

\newtheorem{theorem}{Theorem}

\newtheorem{lemma}{Lemma}

\newenvironment{proof}[1][Proof]{\textbf{#1.} }{\ \rule{0.5em}{0.5em}}

\begin{document}
\title{Drawing dynamical and parameter planes of iterative families and methods\thanks{This research was supported by Ministerio de Ciencia y Tecnolog\'ia MTM2011-28636-C02-02.}
}

\author{Francisco I. Chicharro$^1$, Alicia Cordero$^2$, Juan R. Torregrosa$^3$ \\
\small{$^1$ Instituto de Telecomunicaciones y Aplicaciones Multimedia,}\\
\small{$^{2-3}$ Instituto de Matem\'atica Multidisciplinar}\\
\small{Universitat Polit\`ecnica de Val\`encia,} \\
\small{Camino de Vera, s/n, 46022 Val\`encia, Spain}\\
\small{\{$^1$frachilo@upvnet, $^2$acordero@mat, $^3$jrtorre@mat\}.upv.es}}
\date{}
\maketitle

\begin{abstract}
In this paper  the complex dynamical analysis of the parametric
fourth-order Kim's iterative family is made on quadratic
polynomials, showing the Matlab codes generated to draw the fractal
images necessary to complete the study. The parameter spaces
associated to the free critical points have been analyzed, showing
the stable (and unstable) regions where the selection of the
parameter will provide us excellent schemes (or dreadful ones).\\
\noindent {\bf Keywords}: Nonlinear equation, Kim's family, parameter space, dynamical plane, stability.
\end{abstract}

\section{Introduction}\label{S1}

It is usual to find nonlinear equations in the modelization of many
scientific and engineering problems, and a broadly extended tool to
solve them are the iterative methods. In the last decades, it has
become an increasing and fruitful area of research. More recently,
complex dynamics has revealed as a very useful tool to deep in the
understanding of the rational functions that rise when an iterative
scheme is applied to solve the nonlinear equation $f(z)=0$, with
$f:\mathbb{C}\rightarrow \mathbb{C}$. The dynamical properties of
this rational function give us important information about numerical
features of the method as its stability and reliability.

There is an extensive literature on the study of iteration of
rational mappings of complex variables (see \cite{DH,CGS}, for
instance). The simplest and more deeply analyzed model is the
obtained when $f(z)$ is a quadratic polynomial and the iterative
process is Newton's one. The dynamics of this iterative scheme has
been widely studied (see, among others,
\cite{CGS,Blanchard2,Fallega}).

Just a decade ago that Varona in \cite{V} and Amat et al. in
\cite{Amat3} described the dynamical behavior of several well-known
iterative methods. More recently, in
\cite{GHR,HPR,CCGT,ACCT,SNC,CLND,NSC1,King_AML}, the authors study
the dynamics of different iterative families. In most of these
studies, interesting dynamical planes, including some periodical
behavior and other anomalies, have been obtained. In a few cases,
the parameter planes have been also analyzed.

In order to study the dynamical behavior of an iterative method when
is applied to a polynomial $p(z)$, it is necessary to recall some
basic dynamical concepts. For a more extensive and comprehensive
review of these concepts, see \cite{Blanchard2,devaney}.

Let $R:\mathbb{\hat{C}}\rightarrow\mathbb{\hat{C}}$ be a rational
function, where $\mathbb{\hat{C}}$ is the Riemann sphere. The orbit
of a point $z_0\in\mathbb{\hat{C}}$ is defined as the set of
successive images of $z_0$ by the rational function, $ \left\{
z_0,R(z_0),\ldots,R^n(z_0),\ldots\right\}$.

The dynamical behavior of the orbit of a point on the complex plane
can be classified depending on its asymptotic behavior. In this way,
a point $z_0\in\mathbb{C}$ is a fixed point of $R$ if $R(z_0)=z_0$.
A fixed point is attracting, repelling or neutral if $|R'(z_0)|$ is
lower than, greater than or equal to 1, respectively. Moreover, if
$|R'(z_0)|=0$, the fixed point is superattracting.

If $z^*_f$ is an attracting fixed point of the rational function
$R$, its basin of attraction $\mathcal{A}(z^*_f)$ is defined as the
set of pre-images of any order such that
\begin{equation}
\mathcal{A}(z^*_{f})=\left\{z_0\in\mathbb{\hat{C}}:R^n(z_0)\rightarrow
z^*_{f},n\rightarrow\infty\right\}.\label{cuenca}\nonumber
\end{equation}

The set of points whose orbits tend to an attracting fixed point
$z^*_f$ is defined as the Fatou set, $\mathcal{F}(R)$. The
complementary set, the Julia set $\mathcal{J}(R)$, is the closure of
the set consisting of its repelling fixed points, and establishes
the borders between the basins of attraction.

In this paper, Section \ref{S2} is devoted to the complex analysis
of a known fourth-order family, due to Kim (see \cite{kim}). The
conjugacy classes of its associated fixed point operator, the
stability of the strange fixed points, the analysis of the free
critical points and the analysis of the parameter and dynamical
planes are made. In Section \ref{S4}, the Matlab code used to
generate these tools is showed and the key instructions are
explained in order to help to their eventual modification to adapt
them to other iterative families. Finally, some conclusions are
presented.

\section{Complex dynamics features of Kim's family}\label{S2}
We will focus our attention on the dynamical analysis of a known
parametric family of fourth-order methods for solving a nonlinear
equation $f(x)=0$. Kim in \cite{kim} designs a class of optimal
eighth-order methods, whose two first steps are
\begin{eqnarray}\label{familia}
y_k     & = & x_k - \frac{f(x_k)}{f'(x_k)},\nonumber \\
x_{k+1} & = & y_k - \frac{1+\eta u_k + \lambda u_k^2}{1+(\eta-2)u_k+\mu
u_k^2} \frac{f(y_k)}{f'(x_k)},
\end{eqnarray}
where $u_k=\frac{f(y_k)}{f(x_k)}$. We will suppose $\eta=\mu=0$. The
result is a one-parametric family of iterative schemes whose order
of convergence is four, with no conditions on $\lambda$.

In order to study the affine conjugacy classes of the iterative
methods, the following Scaling Theorem can be easily checked.

\begin{theorem}
Let $g(z)$ be an analytic function, and let $A(z)=\alpha_1 z+\alpha_2$,
with $\alpha_1\neq 0$, be an affine map. Let $h(z)=\gamma(g\circ
A)(z)$, with $ \gamma\neq0$. Let $O_p(z)$ be the fixed point
operator of Kim's family on the polynomial $p(z)$. Then, $(A\circ O_h\circ
A^{-1})(z)=O_g(z)$, that is, $O_g(z)$ and $O_h(z)$ are affine conjugated by
$A$.
\end{theorem}

This result allows up the knowledge of a family of polynomials with
just the analysis of a few cases, from a suitable scaling.

In the following we will analyze the dynamical behavior of the
fourth-order parametric family (\ref{familia}), on the quadratic
polynomial $(z-a)(z-b)$, where $a,b \in \mathbb{C}$.

We apply the M\"{o}bius transformation
$$M(u)=\frac{u-a}{u-b},$$ whose inverse is
$$[M(u)]^{-1}=\frac{u b-a}{u-1},$$ in order to obtain the
one-parametric operator
\begin{equation}\label{op}
O_p(z,\lambda)=-\frac{z^4(1-\lambda+4z+6z^2+4z^3+z^4)}{-1-4z-6z^2-4z^3+(-1+\lambda)z^4},
\end{equation}
associated to the iterative method. In the study of the rational
function (\ref{op}), $z=0$ and $z=\infty$ appear as superattracting
fixed points and $z=1$ is a strange fixed point for $\lambda \neq 1$
and $\lambda \neq 16$. There are also six strange fixed points more (a fixed point is called strange if it does not correspond to any root of the polynomial),
whose analytical expression, depending on $\lambda$, is very
complicated.

As we will see in the following, not only the number but also the
stability of the fixed points depend on the parameter of the family.
The expression of the differentiated operator, needed to analyze the
stability of the fixed points and to define the critical points, is
\begin{equation}
O_p'(z,\lambda)=-\frac{4z^3(1+z)^4(-(1+z)^4+\lambda(1-z+z^2-z^3+z^4))}{(1+4z+6z^2+4z^3-(-1+\lambda)z^4)^2}.
\end{equation}

As they come from the roots of the polynomial, it is clear that the
origin and $\infty $ are always superattractive fixed points, but
the stability of the other fixed points can change depending on the
values of the parameter $\lambda$. In the following result we
establish the stability of the strange fixed point $z=1$.

\begin{theorem}\label{teoestabilidad1}
The character of the strange fixed point  $z=1$ is:
\begin{itemize}
\item[i)] If $\left|\lambda-16\right|>64$, then $z=1$ is an attractor and it cannot be a superattractor.
\item[ii)] When $\left|\lambda-16\right|=64$,  $z=1$ is a parabolic point.
\item[iii)] If $\left|\lambda-16\right|<64$, being $\lambda \neq 1$ and $\lambda \neq 16$, then $z=1$ is a repulsor.
\end{itemize}
\end{theorem}
\proof It is easy to prove that
\[
O_{p}'\left( 1,\lambda\right) =\frac{64}{16-\lambda}.
\]
So,
\[
\left|\frac{64 }{16-\lambda}\right| \leq 1\phantom{aaa} \mbox{is
equivalent to} \phantom{aaa} 64 \leq \left|16-\lambda\right| .
\]
Let us consider $\lambda =a+ib$ an arbitrary complex number. Then,
\[
64^2 \leq 16^2-32a+a^{2}+b^{2}.
\]
That is,
\[
  \left( a-16\right) ^{2}+b^{2}\geq 64^2.
\]
Therefore,
\[
\left| O_{p}^{\prime }\left( 1,\lambda\right) \right| \leq 1 \ \
\mbox{if and only if} \ \ \left|\lambda-16\right|\geq 64.
\]
Finally, if $\lambda$ verifies $\left|\lambda-16\right|< 64$,
then $\left| O_{p}^{\prime }\left( 1,\lambda\right) \right|
>1$ and $z=1$ is a repulsive point, except if $\lambda=1$ or $\lambda=16$, values for which $z=1$ is not a fixed point.
\endproof

The critical points of $O_p(z,\lambda)$ are $z=0$, $z=\infty$, $z=-1$ and
\begin{eqnarray*}
cr_1(\lambda)&=&\frac{1}{4}\left[1+\frac{1}{\lambda-1}-\beta - \sqrt{2} \sqrt{\frac{-5\lambda(6-7\lambda+\lambda^2)}{(\lambda-1)^3}-\frac{(4+\lambda)\beta}{\lambda-1}} \right],\\
cr_2(\lambda)&=&\frac{1}{4}\left[1+\frac{1}{\lambda-1}-\beta + \sqrt{2} \sqrt{\frac{-5\lambda(6-7\lambda+\lambda^2)}{(\lambda-1)^3}-\frac{(4+\lambda)\beta}{\lambda-1}} \right],\\
cr_3(\lambda)&=&\frac{1}{4}\left[1+\frac{1}{\lambda-1}+\beta - \sqrt{2} \sqrt{\frac{-5\lambda(6-7\lambda+\lambda^2)}{(\lambda-1)^3}+\frac{(4+\lambda)\beta}{\lambda-1}} \right],\\
cr_4(\lambda)&=&\frac{1}{4}\left[1+\frac{1}{\lambda-1}-\beta + \sqrt{2} \sqrt{\frac{-5\lambda(6-7\lambda+\lambda^2)}{(\lambda-1)^3}+\frac{(4+\lambda)\beta}{\lambda-1}} \right],\\
\end{eqnarray*}
where $\lambda \neq 1$ and $\beta=\displaystyle\frac{\sqrt{5} \sqrt{(\lambda-1)^2}\sqrt{\lambda (4+\lambda)}}{(\lambda-1)^2}$.

The relevance of the knowledge of the free critical points (critical points different from the roots) is the following known fact: each invariant Fatou component is associated with, at least, one critical point.

\begin{lemma}\label{lemacriticos}
Analyzing the equation $O'_{p}(z, \lambda)=0$, we obtain
\begin{itemize}
\item[a)] If $\lambda=0$, there is no free critical points of operator $O_{p}\left(z,0\right)$.
\item[b)] If $\lambda=16$, then there are four free critical points: $z=-1$, $cr_1\left(16 \right)=\frac{1}{4}\left( -\frac{4}{3}-\frac{8\sqrt{2}}{3}i \right)$, $cr_2\left(16\right)=\frac{1}{4}\left( -\frac{4}{3}+\frac{8\sqrt{2}}{3}i \right)$ and $cr_3\left(16\right)=cr_4\left(16\right)=1$.
\item[c)] If $\lambda=-4$, then there are three different critical points:  $z=-1$, $cr_1\left(-4\right)=cr_3\left(-4\right)=-i$ and $cr_2\left(-4\right)=cr_4\left(-4\right)=i$.
\item[d)] In case of $\lambda=1$, the set of free critical points is $\{ -1,-\frac{1}{2}+\frac{\sqrt{3}}{2}i ,-\frac{1}{2}-\frac{\sqrt{3}}{2}i \}$.
\item[e)] In any other case, $z=-1$,  $cr_1\left(\lambda\right)$, $cr_2\left(\lambda\right)$, $cr_3\left(\lambda\right)$ and $cr_4\left(\lambda\right)$ are the free critical points.
\end{itemize}
Moreover, it can be proved that all free critical points are not independent, as $cr_1(\lambda)=\frac{1}{cr_2(\lambda)}$ and $cr_3(\lambda)=\frac{1}{cr_4(\lambda)}$.
\end{lemma}
Some of these properties determine the complexity of the operator,
as we can see in the following results.
\begin{theorem}
The only member of the family whose operator is always conjugated
to the rational map $z^4$ is the element corresponding to $\lambda=0$.
\end{theorem}
\proof From (\ref{op}), we denote
$p(z)=(1 + z)^4 (-(1 + z)^4 + \lambda (1 - z + z^2 - z^3 + z^4))$ and
$q(z)=(1 + 4 z + 6 z^2 + 4 z^3 - (-1 + \lambda) z^4)^2$. By factorizing both polynomials,
we can observe that the unique value of $\lambda$ verifying
$p(z)=q(z)$ is $\lambda=0$.
\endproof

In fact, the element of Kim's class corresponding to $\lambda=0$ is Ostrowski's method. So, it is the
most stable scheme of the family, as there are no free critical point and the iterations can only
converge to any of the images of the roots of the polynomial. This is the same behavior observed when Ostrowski's scheme was analyzed by the authors as a member of King's family in \cite{King_AML}.

\begin{theorem}\label{teo_orden5}
The element of the family corresponding to $\lambda=1$ is a fifth-order method whose operator is the rational map
\begin{equation}\label{op l=1}
O_p(z,\lambda)=\frac{z^5(2+z)(2+2z+z^2)}{(1+2z)(1+2z+2z^2)}.
\end{equation}
\end{theorem}
\proof From directly substituting $\lambda =1$ in the rational operator (\ref{op}), expression (\ref{op l=1}) is obtained, showing that $z=1$ is not a fixed point in this particular case. Moreover,
\begin{equation*}
O_p'(z,\lambda)=\frac{20 z^4 (1+z)^4 \left(1+z+z^2\right)}{(1+2 z)^2 \left(1+2 z+2 z^2\right)^2},
\end{equation*}
and there exist only three free critical points.
\endproof

Then, in the particular case $\lambda =1$, the order of convergence
is enhanced to five and, although there are three free critical
points, they are in the basin of attraction of zero and infinity,
as the strange fixed points are all repulsive in this case. So, it
is a very stable element of the family with increased convergence in
case of quadratic polynomials.

\subsection{Using the parameter and dynamical planes}

From the previous analysis, it is clear that the dynamical behavior of the rational operator associated to each value of the parameter can be very different. Different parameter spaces associated with
 free critical points of this family are obtained. The process to obtain these parameter planes is the following: we associate each point of the parameter plane with a complex value
of $\lambda$, i.e., with an element of family (\ref{familia}). Every
value of $\lambda$ belonging to the same connected component of the
parameter space gives rise to subsets of schemes of family
(\ref{familia}) with similar dynamical behavior. So, it is interesting
to find regions of the parameter plane as much stable as possible,
because these values of $\lambda$ will give us the best members of the
family in terms of numerical stability.

\begin{figure}[ht!]
\begin{center}
  \includegraphics[width=.75\linewidth]{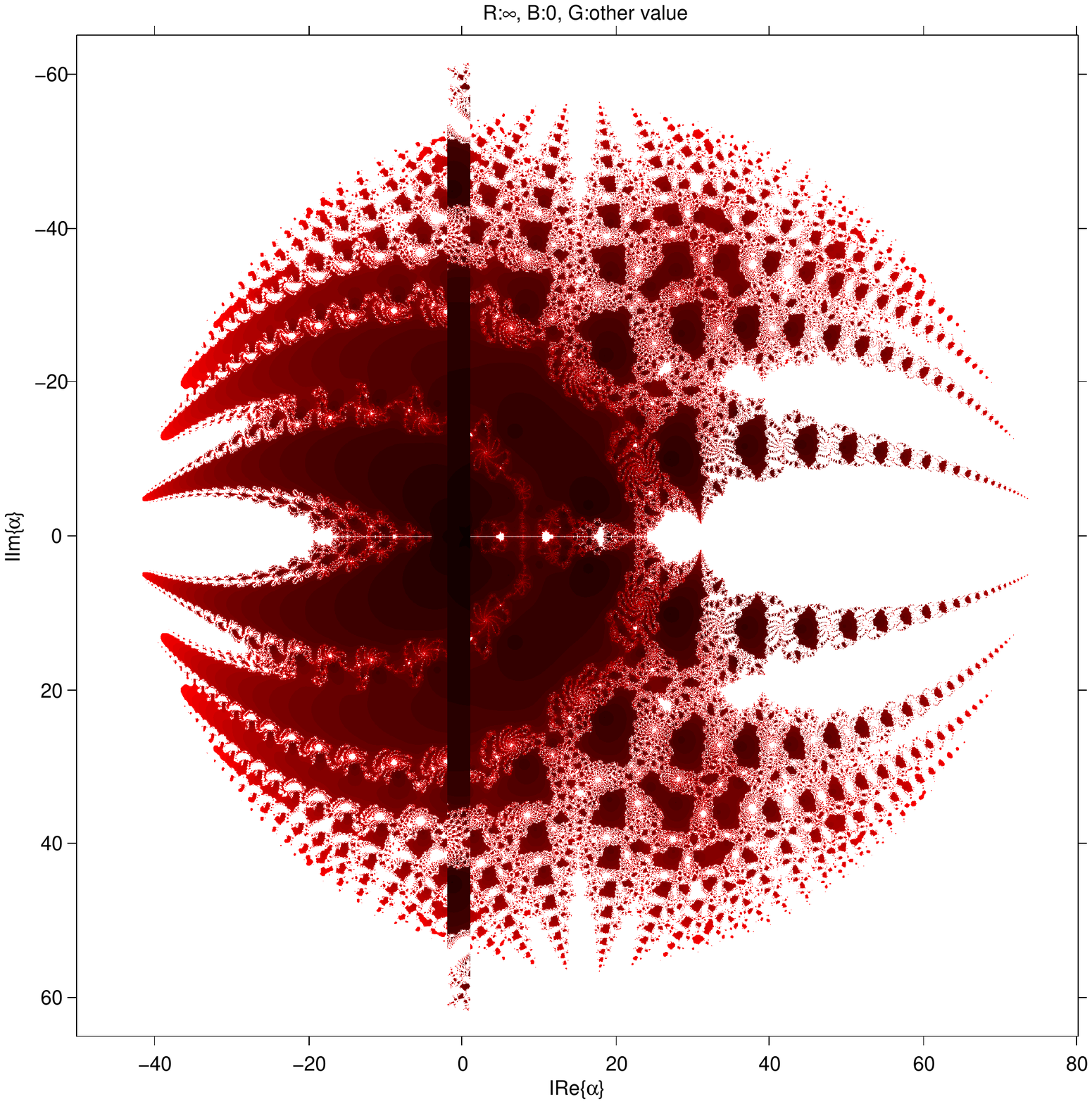}
  \caption{Parameter plane $P_1$ associated to $z=cr_i(\lambda)$, $i=1,2$}\label{planopar1}
  \end{center}
\end{figure}

As $cr_1(\lambda)=\frac{1}{cr_2(\lambda)}$ and $cr_3(\lambda)=\frac{1}{cr_4(\lambda)}$ (see Lemma \ref{lemacriticos}),
we have at most three independent free critical points. Nevertheless, $z=-1$ is pre-image of
the fixed point $z=1$ and the parameter plane associated to this critical point is not significative.
So, we can obtain two different parameter planes, with complementary
information. When we consider the free critical point $cr_1(\lambda)$ (or $cr_2(\lambda)$) as a
starting point of the iterative scheme of the family associated to
each  complex value of $\lambda$, we paint this point of the complex
plane in red if the method converges to any of the roots (zero and
infinity) and they are white in other cases. Then, the parameter
plane $P_1$  is obtained; it is showed in Figure \ref{planopar1}.

\begin{figure}[h!]
  \centering
  \subfloat[$\left|\lambda-16\right|> 64$]{\label{figura5}\includegraphics[width=0.45\textwidth]{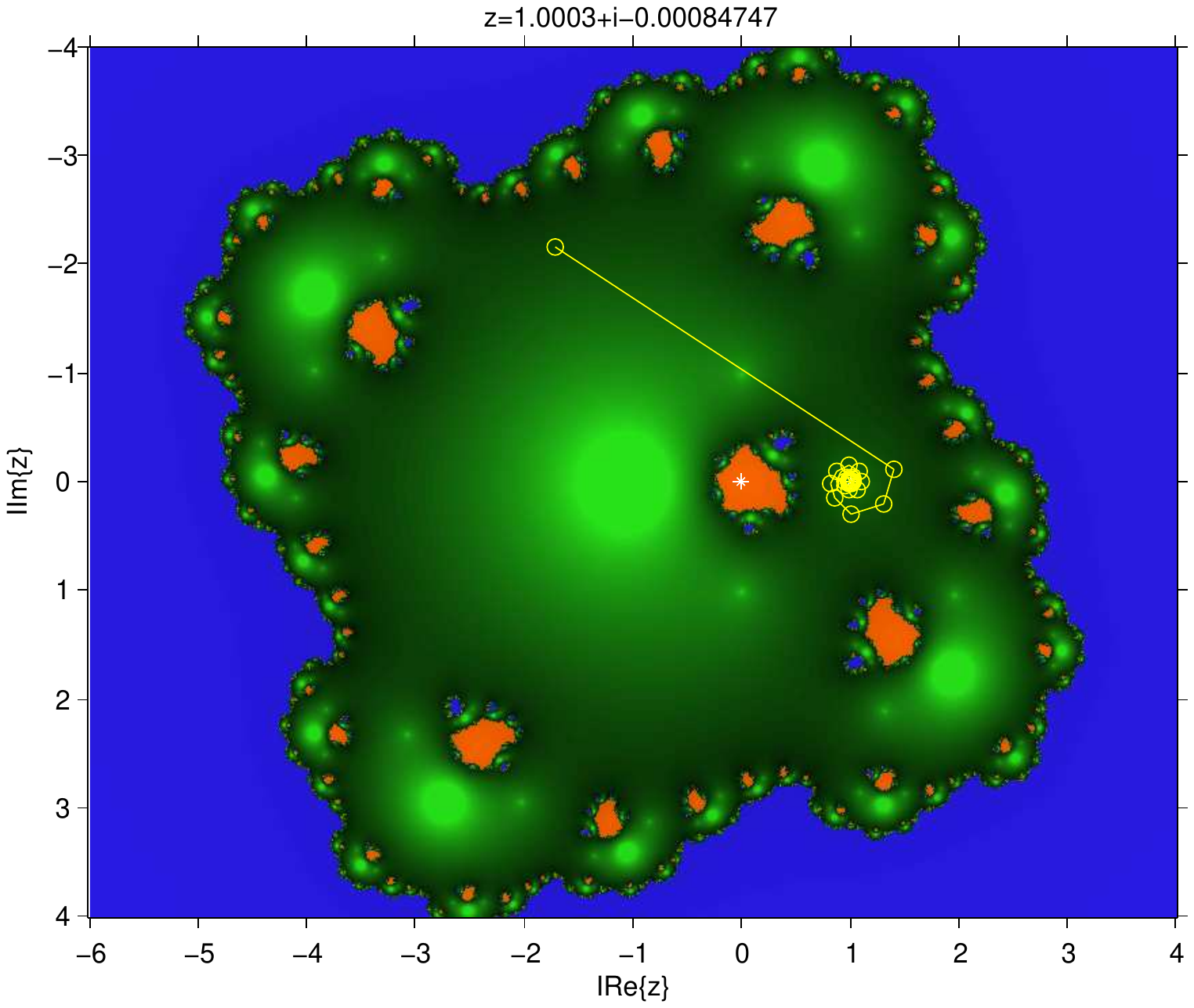}}
  \subfloat[$\left|\lambda-16\right|= 64$]{\label{figura6}\includegraphics[width=0.45\textwidth]{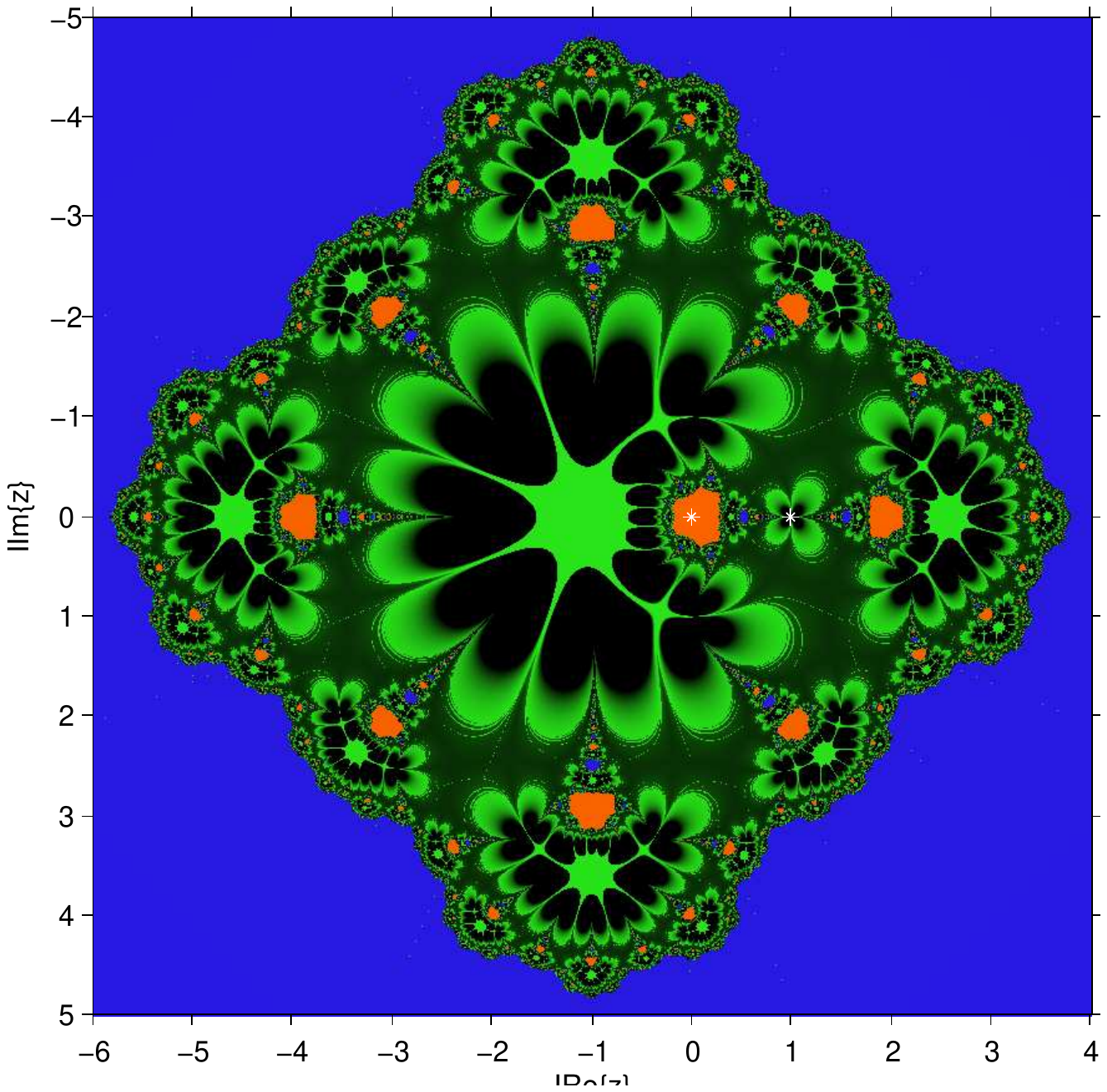}}
    \caption{Dynamical planes for $\lambda$ verifying $\lambda$ verifying $\left|\lambda-16\right|\geq 64$}
\end{figure}

This figure has been generated for values of $\lambda$ in $[-50,80]\times[-65,65]$, with a mesh of $2000\times2000$ points and 400 iterations per point. In $P_1$ the disk of repulsive behavior of $z=1$ is observed, showing different white regions where the convergence to $z\neq0$, $z \neq \infty$ has been reached. An example of a dynamical plane associated to a value of the parameter is shown in Figure \ref{figura5}, where three different basins of attraction appear, two of them of the superattractors $0$ and $\infty$ and the other of $z=1$, that is a fixed attractive point. It can be observed how the orbit (in yellow in the figure) converges asymptotically to the fixed point. Also in Figure \ref{figura6}, the behavior in the boundary of the disk of stability of $z=1$ is presented, where this fixed point is parabolic. An orbit would tend to the parabolic point alternating two "sides" (up and down the parabolic point, in this case).

The generation of dynamical planes is very similar to the one of parameter spaces. In case of dynamical planes, the value of parameter $\lambda$ is constant (so, the dynamical plane is associated to a concrete element of the family of iterative methods). Each point of the complex plane is considered as a starting point of the iterative scheme and it is painted in different colors depending on the point which it has converged to. A detailed explanation of the generation of these graphics, joint with the Matlab codes used to generate them is provided in Section \ref{S4}.

%

\begin{figure}[h!]
  \centering
  \subfloat[Detail of $P_1$]{\label{figura3}\includegraphics[width=0.45\textwidth]{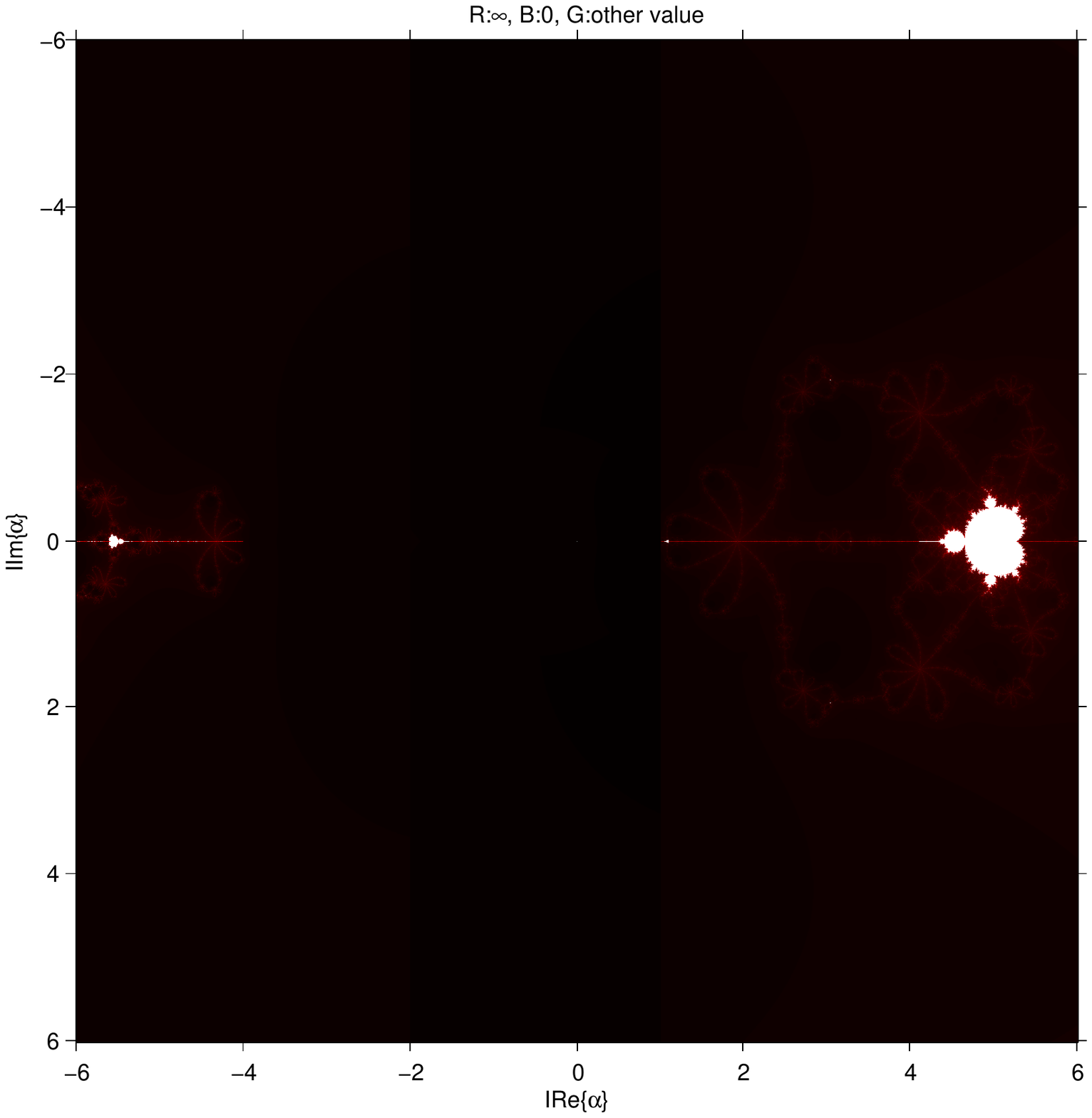}}
  \subfloat[Dynamical plane for $\lambda=1$]{\label{figura4}\includegraphics[width=0.45\textwidth]{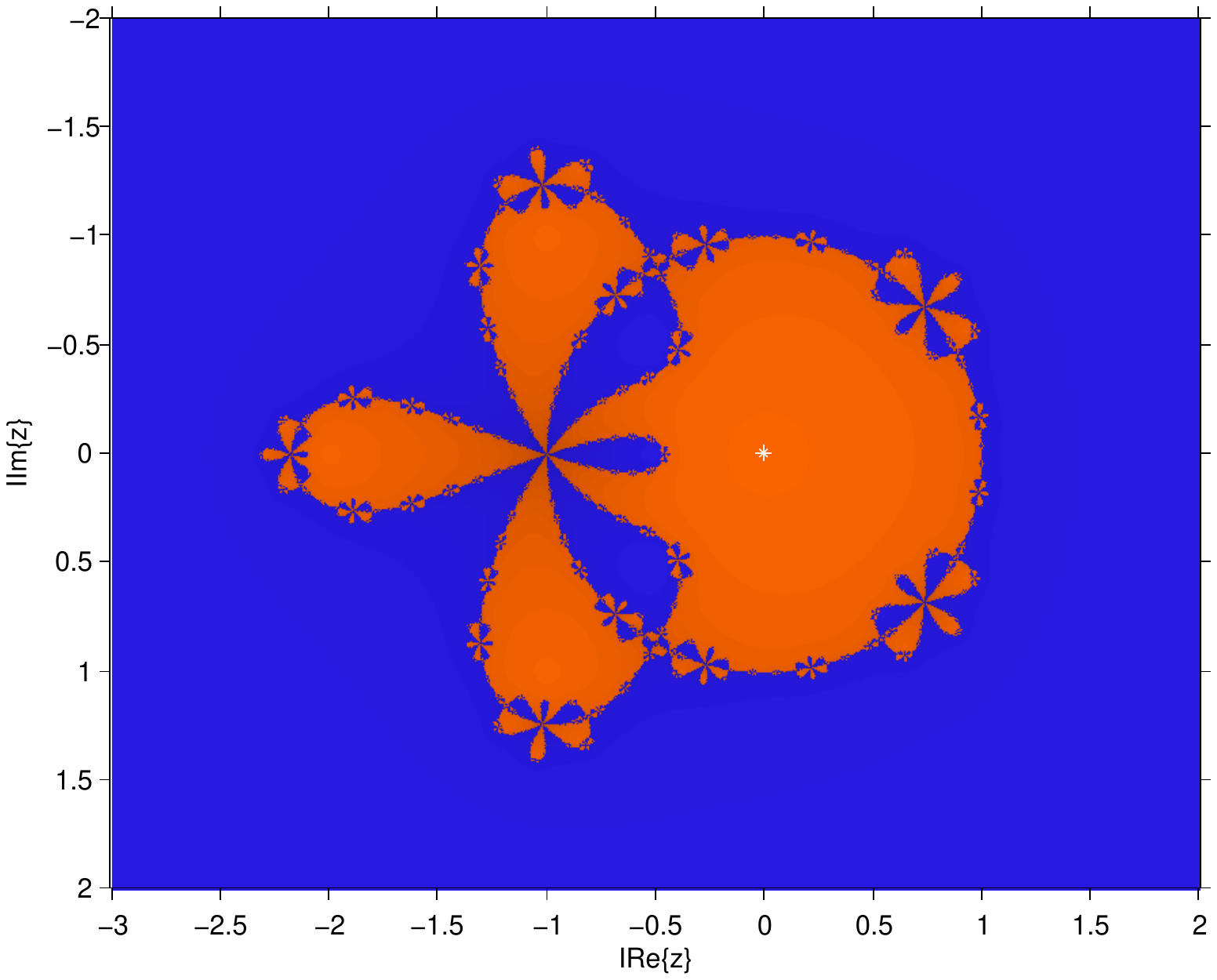}}
    \caption{Around the origin}
\end{figure}

In Figure \ref{figura3}, a detail of the region around $\lambda=0$ is seen. Let us
notice that region around the origin is specially stable, specifically the vertical band between $-4$ and $1$ (see also Figure \ref{figura6}).

%
%

 In fact, for $\lambda=0$, the associated dynamical plane is just a disk and its complementary in $\mathbb{C}$, as in Newton's method. Around the origin is also very stable, with two connected components in the Fatou set. When $\lambda=16$, $z=1$ is not a fixed point (see Theorem \ref{teoestabilidad1}) and $\{-1,1\}$ defines a periodic orbit of period 2 (see Figure \ref{figura4}). The singularity of this value of the parameter can be also observed in Figure \ref{figure7}(AÑADIR), in which a dynamical plane for $\lambda=15.9-0.2i$ is presented, showing a very stable behavior with only two basins of attraction, corresponding to the image of the roots of the polynomial by the M\"{o}bius map.
 \begin{figure}[h!]
\begin{center}
  \includegraphics[width=.6\linewidth]{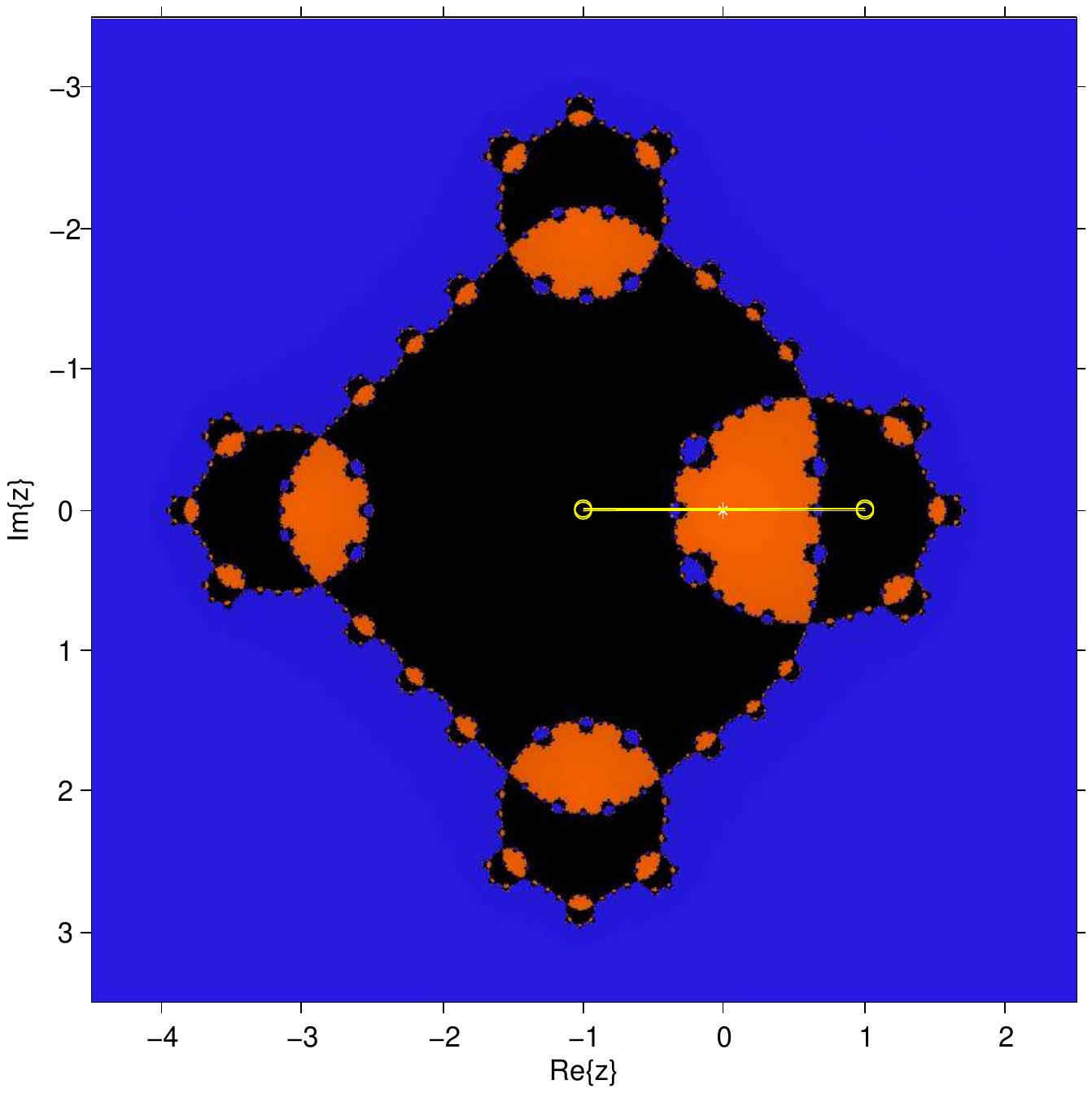}
  \caption{2-periodic orbit for $\lambda=16$}\label{figura4}
   \end{center}
\end{figure}

It is also interesting to note in Figure \ref{figura3} that white figures with a certain similarity with
the known Mandelbrot set appear. Their antennas end in the values $\lambda=-4$ and $\lambda=1$,
whose dynamical behavior is very different from the near values of the parameter, as was shown in Lemma \ref{lemacriticos}.

\begin{figure}[h!]
\begin{center}
  \includegraphics[width=.75\linewidth]{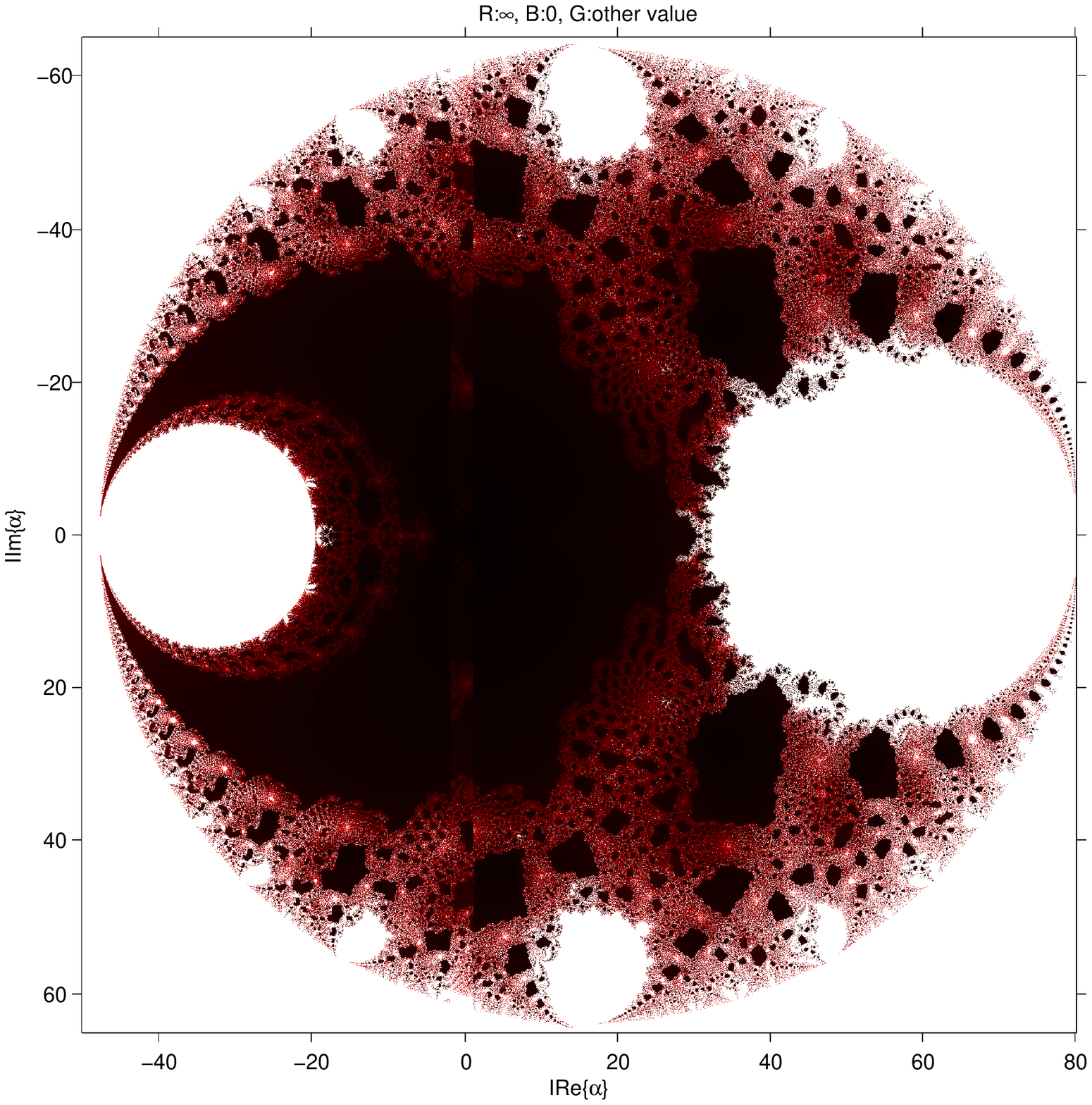}\\
  \caption{Parameter space $P_2$, associated to $z=cr_i$, $i=3,4$}\label{planopar2}
   \end{center}
\end{figure}

A similar procedure can be carried out with the free critical
points, $z=cr_i$, $i=3,4$, obtaining the parameter planes $P_2$,
showed in Figure \ref{planopar2}.

\begin{figure}[h!]
  \centering
  \subfloat[2-periodic orbit]{\label{figura7}\includegraphics[width=0.45\textwidth]{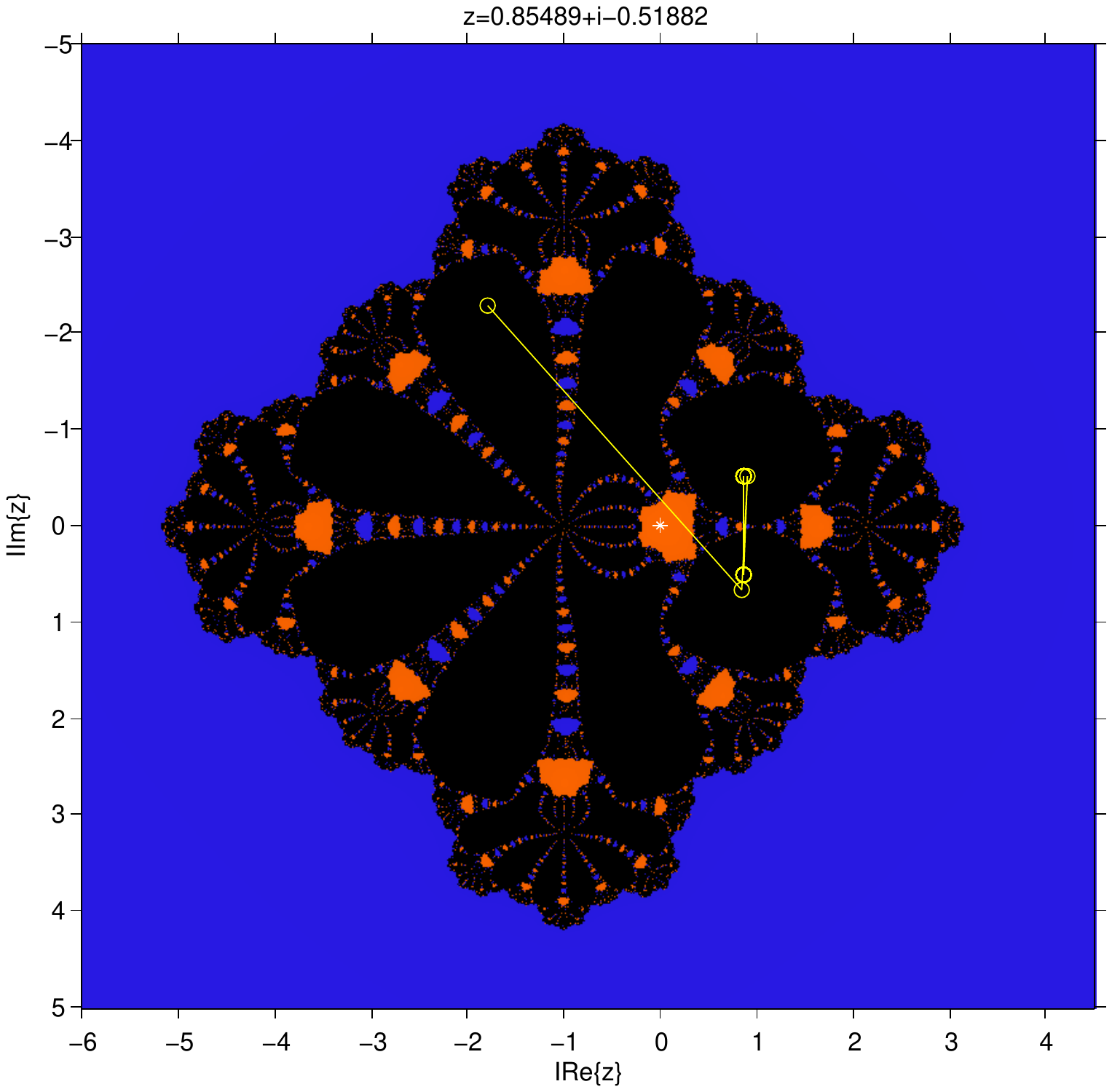}}
  \subfloat[Two atracting strange fixed points]{\label{figura8}\includegraphics[width=0.45\textwidth]{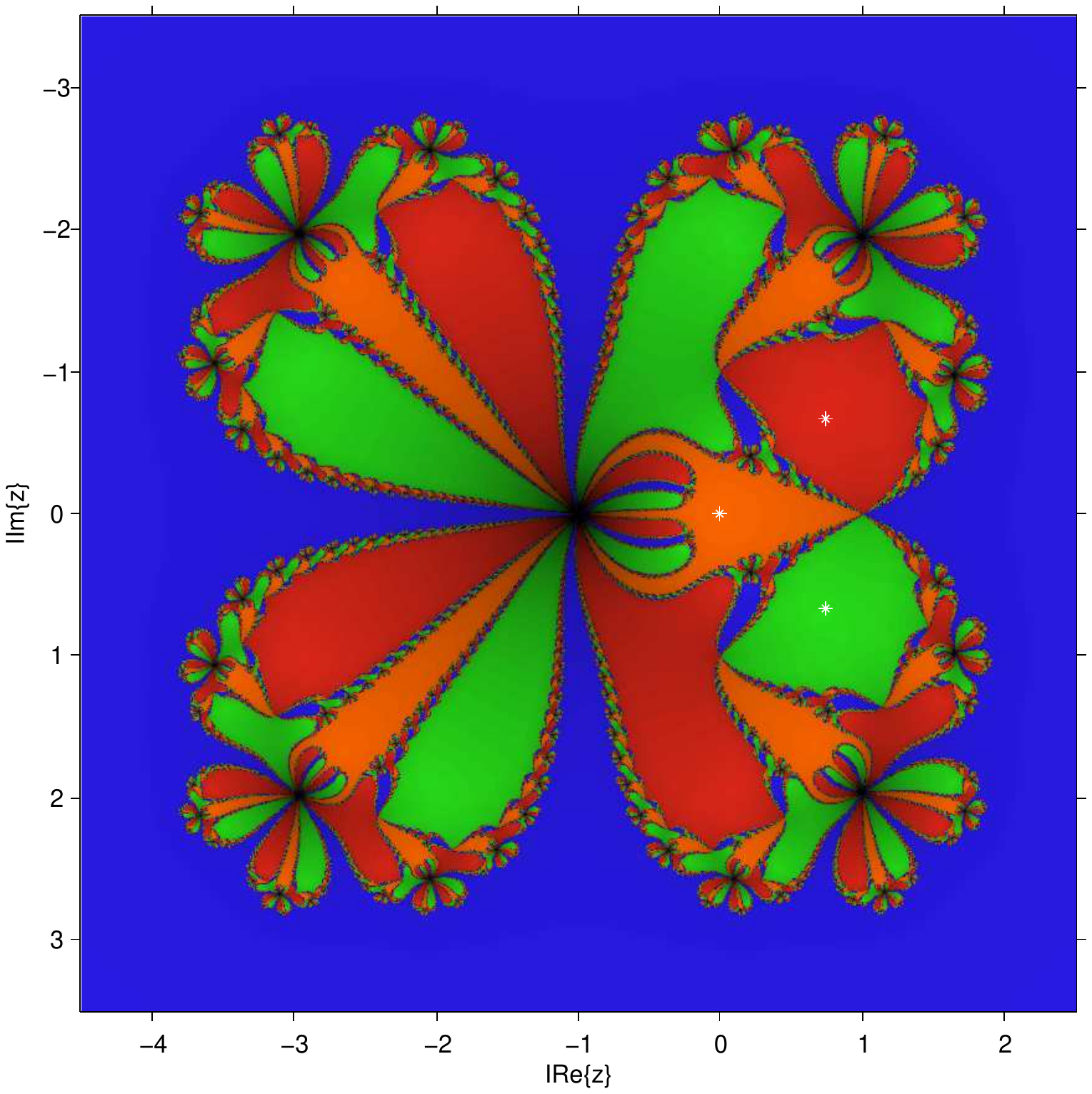}}\\
   \subfloat[4-periodic orbit]{\label{figura9}\includegraphics[width=0.45\textwidth]{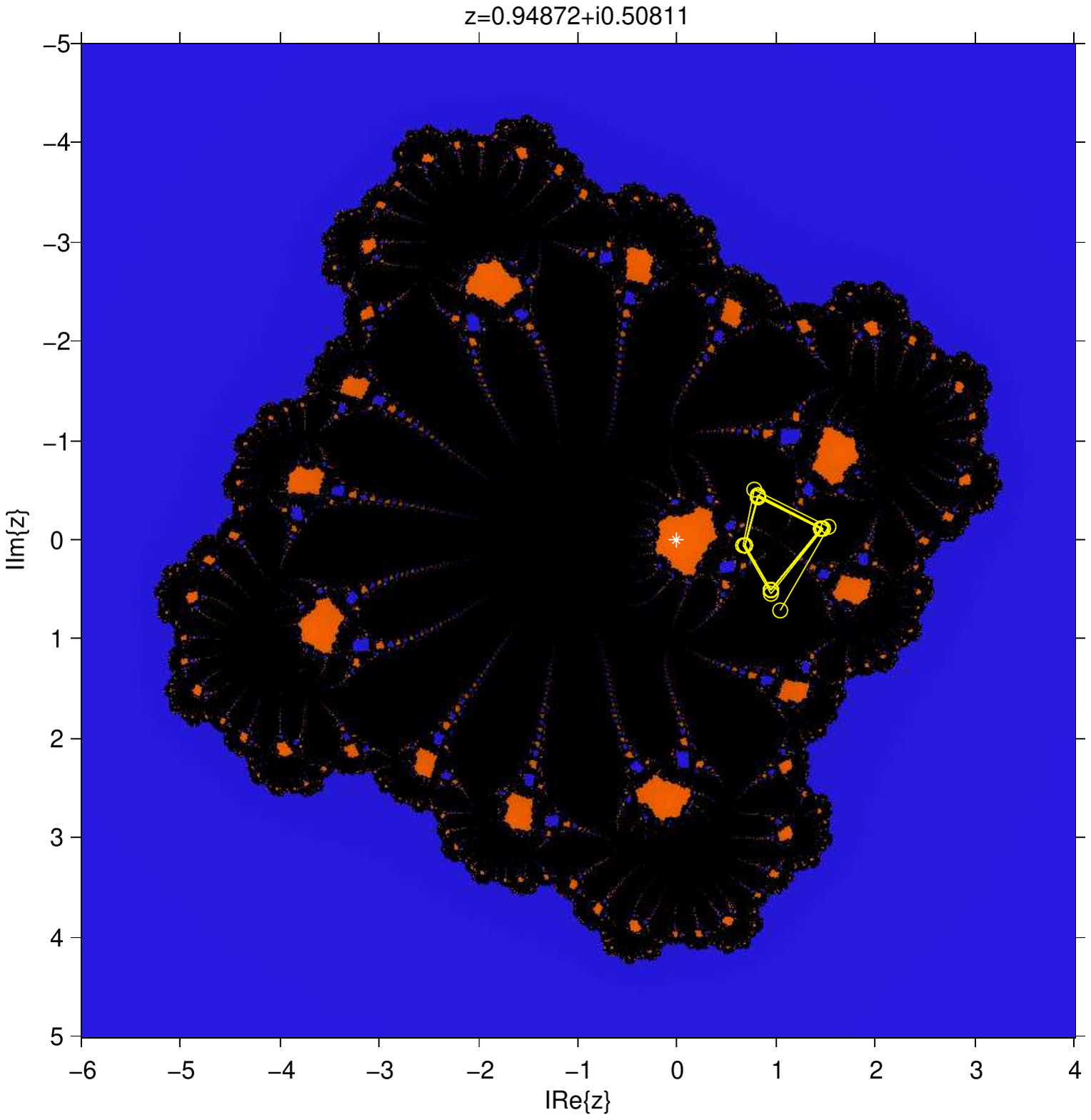}}
  \subfloat[3-periodic orbit]{\label{figura10}\includegraphics[width=0.45\textwidth]{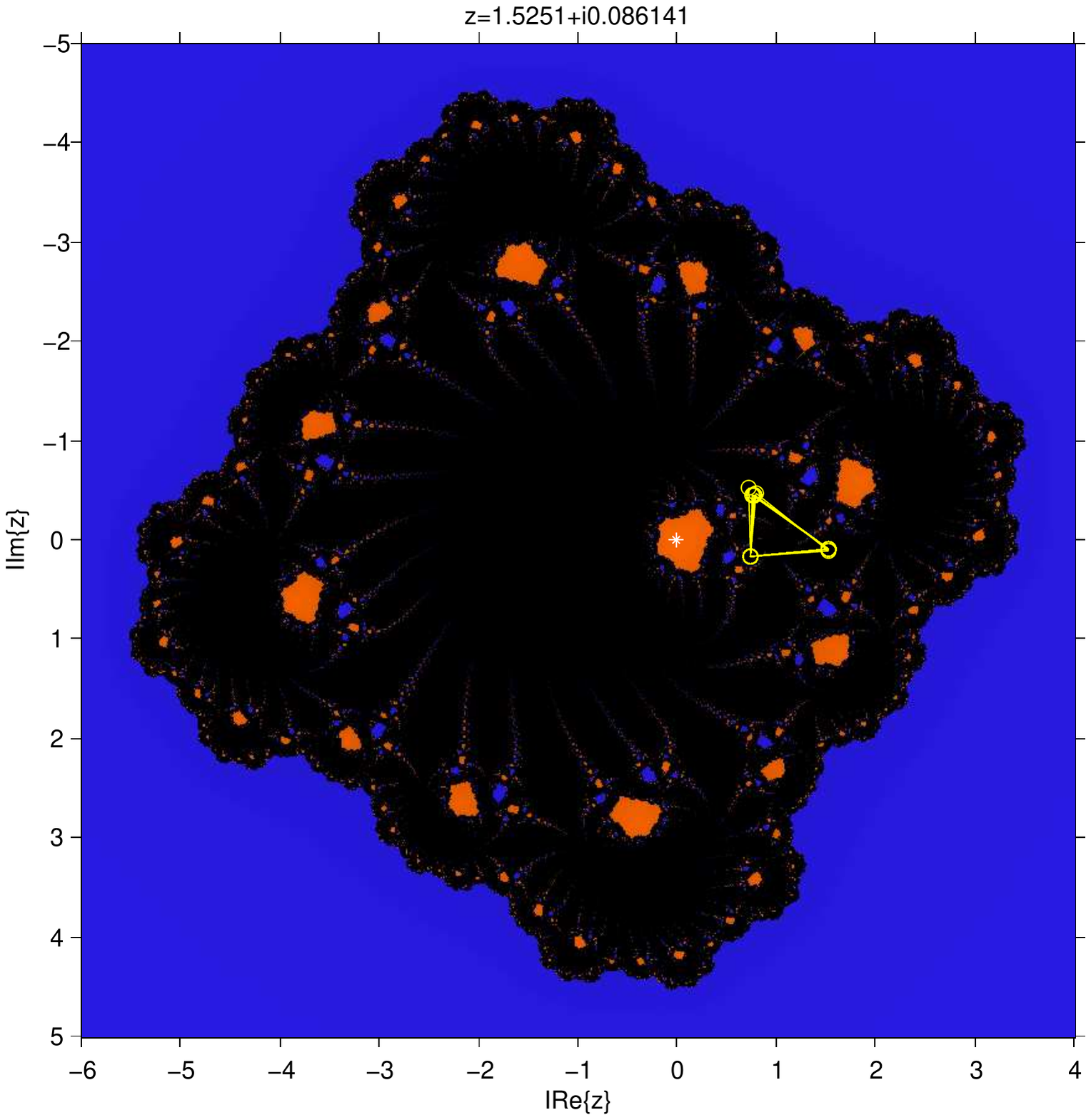}}
    \caption{Some dynamical planes from $P_2$ }
\end{figure}

 As in case of $P_1$, the disk of repulsive behavior of $z=1$ is clear, and inside it different "bulbs" appear, similar to disks. The biggest on the left of the real axis corresponds to the set of values of $\lambda$ where the fixed point $z=1$ has bifurcated in a periodic orbit of period two, as can be seen in Figure \ref{figura7}. In the right of the real axis a bulb is the loci of two conjugated strange fixed points, see Figure \ref{figura8}.

The bulbs on the top (see Figure \ref{figura9}) and on the bottom of the imaginary axis correspond to periodic orbits of period 4. The rest of the bulbs surrounding the boundary of the stability disk of $z=1$ correspond to regions where periodic orbits of different periods appear. In fact, we can observe in  Figure \ref{figura10}) a periodic orbit of period 3, obtained from $\lambda=50+50i$ (comprobar). By applying Sharkovsky's Theorem (see \cite{devaney}), we can affirm that periodic orbits of arbitrary period can be found.

\section{MATLAB$^{\copyright}$ planes code}\label{S4}
The main goal of drawing the dynamical and parameters planes is the comprehension of the family or method behavior at a glance. The procedure to generate a dynamical or a parameter plane is very similar. However, there are significative differences, so both cases are developed below.

\subsection{Dynamical planes}
From a fixed point operator -- that associates a polynomial to an iterative method -- the dynamical plane illustrates the basins of attraction of the operator. The orbit of every point in the dynamical plane tends to a root (or to the infinity); this information, and the speed with which that the points tends to the root, can be displayed in the dynamical plane.
In our pictures, each basin of attraction is drawn with a different color. Moreover, the brightness of the color points the number of iterations needed to reach the root of the polynomial.

The following code covers the Kim's fixed point operators, when it is applied to a quadratic polynomial. This code has been utilized to generate the dynamical planes of several papers, as \cite{CCGT}, \cite{ACCT} or \cite{King_AML}. (poner gato)
\newline

\DefineShortVerb[fontfamily=courier,fontseries=m]{\$}
\DefineShortVerb[fontfamily=courier,fontseries=b]{\#}

\noindent\scriptsize
 \hspace*{-1.6em}{\tiny 1}$ $ \color{mblue}$function$ \color{black}$ [I,it]=dynamicalPlane(lambda,bounds,points,maxiter)$\\
 \hspace*{-1.6em}{\tiny 2}$  $\\
  \hspace*{-1.6em}{\tiny 3}$  $\color{mgreen}$ \%\% Description $\color{black}\\
  \hspace*{-1.6em}{\tiny 4}$  $\color{mgreen}$ \% - dynamicalPlane obtains the dynamical plane of the Kim iterative method $\color{black}\\
 \hspace*{-1.6em}{\tiny 5}$   $\color{mgreen}$ \% when it is applied to a quadratic polynomial. The dynamical plane is obtained $\color{black}\\
 \hspace*{-1.6em}{\tiny 6}$   $\color{mgreen}$ \% as a 'points'-by-'points'-by-3 matrix 'I', and can be displayed as $\color{black}\\
 \hspace*{-1.6em}{\tiny 7}$   $\color{mgreen}$ \% >> imshow(I); $\color{black}\\
 \hspace*{-1.6em}{\tiny 8}$   $\color{mgreen}$ \% moreover, the 'points'-by-'points matrix 'it' records the number of $\color{black}\\
 \hspace*{-1.6em}{\tiny 9}$   $\color{mgreen}$ \% iterations of each point. $\color{black}\\
 \hspace*{-1.6em}{\tiny 10}$   $\color{mgreen}$ \% - the method is iterated till the 'maxiter' iterations is reached, or $\color{black}\\
 \hspace*{-1.6em}{\tiny 11}$   $\color{mgreen}$ \% till the estimation is enough close to the root. $\color{black}\\
 \hspace*{-1.6em}{\tiny 12}$   $\color{mgreen}$ \% - it is mandatory the previous execution of $\color{black}\\
 \hspace*{-1.6em}{\tiny 13}$   $\color{mgreen}$ \%   >> syms x $\color{black}\\
 \hspace*{-1.6em}{\tiny 14}$   $\color{mgreen}$ \% bounds:   [min(Re(z)) max(Re(z)) min(Im(z)) max(Im(z))] $\color{black}\\
 \hspace*{-1.6em}{\tiny 15}$   $\color{mgreen}$ \% test:     [I,it]=dynamicalPlane(0,[-1 1 -1 1],400,20); $\color{black}\\
 \hspace*{-1.6em}{\tiny 16}$   $\color{mgreen}$   \% Values $\color{black}\\
 \hspace*{-1.6em}{\tiny 17}$  x0=bounds(1); xN=bounds(2); y0=bounds(3); yN=bounds(4); $\\
 \hspace*{-1.6em}{\tiny 18}$  funfun=matlabFunction(fun);$\\
 \hspace*{-1.6em}{\tiny 19}$  $\\
 \hspace*{-1.6em}{\tiny 20}$  $\color{mgreen}$\% Fixed Point Operator $\color{black}\\
 \hspace*{-1.6em}{\tiny 21}$  syms $\color{mred}$x z$\color{black}\\
 \hspace*{-1.6em}{\tiny 22}$  $\color{mgreen}$\% Kim's operator $\color{black}\\
 \hspace*{-2em}{\tiny 23}$  Op=simple(-x.^4*(1-lambda+4*x+6*x^2+4*x^3+x^4)/(-1-4*x-6*x^2-4*x^3+x^4*(lambda-1)));$
 \hspace*{-2em}{\tiny 24}$      $\\
 \hspace*{-2em}{\tiny 25}$  $\color{mgreen}$\% Attracting points$\color{black}\\
 \hspace*{-2em}{\tiny 26}$  fOp=matlabFunction(Op);$\\
 \hspace*{-2em}{\tiny 27}$  Opx=Op-x;$\\
 \hspace*{-2em}{\tiny 28}$  pf=double(solve(Opx));$\\
 \hspace*{-2em}{\tiny 29}$  dOp=diff(Op);$\\
 \hspace*{-2em}{\tiny 30}$  pc=double(solve(factor(dOp)));$\\
\hspace*{-2em}{\tiny 31}$   adOp=matlabFunction(abs(dOp));$\\
\hspace*{-2em}{\tiny 32}$   inda=double(abs(adOp(pf)))<=1;$\\
\hspace*{-2em}{\tiny 33}$   pa=pf(inda==1);$\\
\hspace*{-2em}{\tiny 34}$  $ \color{mblue}$if$\color{black}$ isempty(pa)$\\
\hspace*{-2em}{\tiny 35}$     pa=double(solve(fun));$\\
\hspace*{-2em}{\tiny 36}$  $ \color{mblue}$end$\color{black}\\
 \hspace*{-2em}{\tiny 37}$   $\\
 \hspace*{-2em}{\tiny 38}$  $\color{mgreen}$\% Preparing the image$\color{black}\\
 \hspace*{-2em}{\tiny 39}$      $\color{mgreen}$\% The image must have an odd number of points$\color{black}\\
 \hspace*{-2em}{\tiny 40}$      $\color{mblue}$if$\color{black}$(mod(points,2)==0)$\\
 \hspace*{-2em}{\tiny 41}$          points=points+1;$\\
 \hspace*{-2em}{\tiny 42}$      $\color{mblue}$end$\color{black}\\
 \hspace*{-2em}{\tiny 43}$  $\\
 \hspace*{-2em}{\tiny 44}$      $\color{mgreen}$\% Complex mesh of points$\color{black}\\
 \hspace*{-2em}{\tiny 45}$      dx=xN-x0; dy=yN-y0; d=max(dx,dy);$\\
 \hspace*{-2em}{\tiny 46}$      step=d/points;$\\
 \hspace*{-2em}{\tiny 47}$      x=x0:step:xN;$\\
 \hspace*{-2em}{\tiny 48}$      y=y0:step:yN;$\\
 \hspace*{-2em}{\tiny 49}$      [X,Y]=meshgrid(x,y);$\\
 \hspace*{-2em}{\tiny 50}$      z=complex(X,Y);$\\
 \hspace*{-2em}{\tiny 51}$  $\\
 \hspace*{-2em}{\tiny 52}$      $\color{mgreen}$\% Matrix startup$\color{black}\\
 \hspace*{-2em}{\tiny 53}$      it=zeros(size(z));$\\
 \hspace*{-2em}{\tiny 54}$      r1=zeros(size(z)); r2=zeros(size(z)); r3=zeros(size(z));$\\
 \hspace*{-2em}{\tiny 55}$      R=zeros(size(z)); G=zeros(size(z)); B=zeros(size(z));$\\
 \hspace*{-2em}{\tiny 56}$  $\\
 \hspace*{-2em}{\tiny 57}$      $\color{mgreen}$\% Colour of each point$\color{black}\\
 \hspace*{-2em}{\tiny 58}$      [f,col]=size(z);$\\
 \hspace*{-2em}{\tiny 59}$      $\color{mblue}$for$\color{black}$ j=1:f$\\
 \hspace*{-2em}{\tiny 60}$          $\color{mblue}$for$\color{black}$ k=1:col$\\
 \hspace*{-2em}{\tiny 61}$              s=z(j,k); rootfound=0;$\color{black}$$\\
 \hspace*{-2em}{\tiny 62}$              $\color{mblue}$while$\color{black}$ (rootfound==0 && it(j,k)<maxiter)$\\
 \hspace*{-2em}{\tiny 63}$                  s=fOp(s);$\\
 \hspace*{-2em}{\tiny 64}$                  it(j,k)=it(j,k)+1;$\\
 \hspace*{-2em}{\tiny 65}$                  $\color{mblue}$if$\color{black}$ norm([real(s)-real(pa(1)) imag(s)-imag(pa(1))])<1e-3$\\
 \hspace*{-2em}{\tiny 66}$                      r1(j,k)=maxiter-1.5*it(j,k);$\\
 \hspace*{-2em}{\tiny 67}$                      R(j,k)=r1(j,k)/maxiter;$\\
 \hspace*{-2em}{\tiny 68}$                      G(j,k)=r1(j,k)/maxiter*102/255;$\\
 \hspace*{-2em}{\tiny 69}$                      rootfound=1;$\\
 \hspace*{-2em}{\tiny 70}$                $\color{mblue}$else$\color{black}\color{mblue}$if$\color{black}$length(pa)>1&&norm([real(s)-real(pa(2))imag(s)-imag(pa(2))])<1e-3$\\
 \hspace*{-2em}{\tiny 71}$                          r2(j,k)=maxiter-2*it(j,k);$\\
 \hspace*{-2em}{\tiny 72}$                          R(j,k)=r2(j,k)/maxiter*40/255;$\\
 \hspace*{-2em}{\tiny 73}$                          G(j,k)=r2(j,k)/maxiter*80/255;$\\
 \hspace*{-2em}{\tiny 74}$                          B(j,k)=r2(j,k)/maxiter;$\\
 \hspace*{-2em}{\tiny 75}$                          rootfound=1;$\\
 \hspace*{-2em}{\tiny 76}$                      $\color{mblue}$else$\color{black}\color{mblue}$if$\color{black}$ length(pa)>2&&norm([real(s)-real(pa(3)) imag(s)-imag(pa(3))])<1e-3$\\
 \hspace*{-2em}{\tiny 77}$                              r3(j,k)=maxiter-it(j,k);$\\
 \hspace*{-2em}{\tiny 78}$                              R(j,k)=r3(j,k)/maxiter*41/255;$\\
 \hspace*{-2em}{\tiny 79}$                              G(j,k)=r3(j,k)/maxiter*230/255;$\\
 \hspace*{-2em}{\tiny 80}$                              B(j,k)=r3(j,k)/maxiter*56/255;$\\
 \hspace*{-2em}{\tiny 81}$                              rootfound=1;$\\
 \hspace*{-2em}{\tiny 82}$                          $\color{mblue}$end$\color{black}\\
 \hspace*{-2em}{\tiny 83}$                      $\color{mblue}$end$\color{black}\\
 \hspace*{-2em}{\tiny 84}$                  $\color{mblue}$end$\color{black}\\
 \hspace*{-2em}{\tiny 85}$              $\color{mblue}$end$\color{black}\\
 \hspace*{-2em}{\tiny 86}$          $\color{mblue}$end$\color{black}\\
 \hspace*{-2em}{\tiny 87}$      $\color{mblue}$end$\color{black}\\
 \hspace*{-2em}{\tiny 88}$  $\\
 \hspace*{-2em}{\tiny 89}$  $\color{mgreen}$\% Image display$\color{black}\\
 \hspace*{-2em}{\tiny 90}$  I(:,:,1)=R(:,:); I(:,:,2)=G(:,:); I(:,:,3)=B(:,:);$\color{black}\\
 \hspace*{-2em}{\tiny 91}$  figure, imshow(I,$\color{mred}$'Xdata'$\color{black}$,[x0 xN], $\color{mred}$'Ydata'$\color{black}$, [y0 yN])$\\
 \hspace*{-2em}{\tiny 92}$  axis $\color{mred}$on$\color{black}$, axis $\color{mred}$xy$\color{black}$, hold $\color{mred}$on$\color{black}\\
 \hspace*{-2em}{\tiny 93}$  plot(real(pa),imag(pa),$\color{mred}$'w*'$\color{black}$)$\\
 \hspace*{-2em}{\tiny 94}$  xlabel($\color{mred}$'Re\{z\}'$\color{black}$); ylabel($\color{mred}$'Im\{z\}'$\color{black}$);$\\
 \hspace*{-2em}{\tiny 95} $axis xy$

\UndefineShortVerb{\$}
\UndefineShortVerb{\#}

\normalsize
The code is divided into five different parts:

\begin{enumerate}
  \item Values (lines 17--18).\\ The bounds are renamed and the symbolic function introduced as {\ttfamily fun} is translated to an anonymous function, recallable by the output handle.
  \item Fixed point operators (lines 23).\\
  \item Calculation of attractive fixed points (lines 26--36).
  \item Image creation (lines 39--94).\\ Once the fixed point operator and the attracting points are set, the next step consists of the determination of the basins of attraction. The combination of the input parameters {\ttfamily bounds} and {\ttfamily points} set the resolution of the image, and it establishes the mesh of complex points (lines 39--50).

      Lines 58--87 are devoted to assign a color to each starting point. It depends on the basin of attraction and the number of iterations needed to reach the root. If the orbit tends to the attracting point set in the first index of line 35, the point is pictured in orange, as lines 67--69 show; for cases second and third, the point is pictured in blue (lines 72--74) and green (lines 78--80), respectively. Otherwise, the point is not modified, so its color is black.

      As the number of iterations needed to reach convergence increases, its corresponding color gets closer to white (black in the decreasing case). A coefficient in each case (lines 66, 71 and 77) is high if the number of iterations is low, and the RGB values are greater than in the slow orbit instance.
  \item Image display (lines 90--94).\\ The image display is based on the {\ttfamily imshow} command. Images are usually displayed in matrix form (from top to bottom and from left to right). In this case, the image is composed by complex points, so the natural display is the cartesian one (from bottom to top and from left to right). With this purpose, {\ttfamily axis xy} is written in line 95.
\end{enumerate}

Once the program is executed, the output values are the image {\ttfamily I} and the number of iterations of each point {\ttfamily it}. Our recommendation is the use of the {\ttfamily surf} command to plot the number of iterations, in combination with the {\ttfamily shading} one.

In order to apply the introduced code to different fixed point operators, the only part to be changed is the fixed point operators corresponding one. If the method can converge to more than three points, just add another {\ttfamily else if} structure (as lines 79--84) and set a color as many times as necessary.

\subsection{Parameter planes}


\DefineShortVerb[fontfamily=courier,fontseries=m]{\$}
\DefineShortVerb[fontfamily=courier,fontseries=b]{\#}
\noindent\scriptsize
 \hspace*{-2em}{\tiny }$ $ \color{mblue}$function$ \color{black}$ [a,I,c]=parametricplane(axini,axfin,ayini,ayfin,points,maxiter)$\\
 \hspace*{-1.6em}{\tiny 2}$  $\\
  \hspace*{-1.6em}{\tiny 3}$  $\color{mgreen}$ \%\% Description $\color{black}\\
  \hspace*{-1.6em}{\tiny 4}$  $\color{mgreen}$ \% - parametricplane obtains the parametric plane of the Kim iterative family $\color{black}\\
 \hspace*{-1.6em}{\tiny 5}$   $\color{mgreen}$ \% when it is applied to a quadratic polynomial, associated with the free critical point cr_2.$\color{black}\\
 \hspace*{-1.6em}{\tiny 6}$   $\color{mgreen}$ \% [axini,axfin,ayini,ayfin] define the rectangle for possible values of the parameter \lambda $\color{black}\\
 \hspace*{-1.6em}{\tiny 7}$   $\color{mgreen}$ \% points defines the mesh of size 'points'-by-'points'$\color{black}\\
 \hspace*{-1.6em}{\tiny 8}$   $\color{mgreen}$ \% maxiter is the maximum number of iterations of the method per value of \lambda $\color{black}\\
 \hspace*{-1.6em}{\tiny 9}$  $ \\
 \hspace*{-1.6em}{\tiny 10}$   $\color{mgreen}$ \% test:     [a,I,c]=parametricplane(-2,2,-2,2,500,25); $\color{black}\\
 \hspace*{-1.6em}{\tiny 11}$   $\color{mgreen}$   \% Values $\color{black}\\
 \hspace*{-2em}{\tiny 12}$   $\\
 \hspace*{-2em}{\tiny 13}$  $\color{mgreen}$\% Preparing the image$\color{black}\\
 \hspace*{-2em}{\tiny 14}$      $\color{mgreen}$\% The image must have an odd number of points$\color{black}\\
 \hspace*{-2em}{\tiny 15}$      $\color{mblue}$if$\color{black}$(mod(points,2)==0)$\\
 \hspace*{-2em}{\tiny 16}$          points=points+1;$\\
 \hspace*{-2em}{\tiny 17}$      $\color{mblue}$end$\color{black}\\
 \hspace*{-2em}{\tiny 18}$  $\\
 \hspace*{-2em}{\tiny 19}$      $\color{mgreen}$\% Complex mesh of points$\color{black}\\
 \hspace*{-1.6em}{\tiny 20}$   ax=linspace(axini,axfin,puntos);$\\
 \hspace*{-1.6em}{\tiny 21}$   ay=linspace(ayini,ayfin,puntos);$\\
 \hspace*{-1.6em}{\tiny 22}$   [AX,AY]=meshgrid(ax,ay);$\\
 \hspace*{-1.6em}{\tiny 23}$   a=complex(AX,AY);$\\
 \hspace*{-2em}{\tiny 24}$  $\\
 \hspace*{-2em}{\tiny 25}$      $\color{mgreen}$\% Matrices startup$\color{black}\\
 \hspace*{-1.6em}{\tiny 26}$   I=zeros(puntos); c=zeros(puntos);$\\
 \hspace*{-1.6em}{\tiny 27}$ R=zeros(puntos); G=zeros(puntos); B=zeros(puntos);$\\
 \hspace*{-2em}{\tiny 57}$      $\color{mgreen}$\% Colour of each point$\color{black}\\
 \hspace*{-2em}{\tiny 30}$      $\color{mblue}$for$\color{black}$ j=1:points$\\
 \hspace*{-2em}{\tiny 31}$          $\color{mblue}$for$\color{black} $ k=1:points$\\
  \hspace*{-2em}{\tiny 32}$                 it=0;$\\
 \hspace*{-2em}{\tiny 33}$                  aa=a(j,k);$\\
 \hspace*{-2em}{\tiny 34}$                  c1=-((-4-aa)/(4*(-1+aa)))-(sqrt(5)*sqrt(4*aa+aa^2))/...$\\
 \hspace*{-2em}{\tiny 35}$                  (4*sqrt(1-2*aa+aa^2))+1/2*sqrt((-4-aa)^2/(2*(-1+aa)^2)...$\\
 \hspace*{-2em}{\tiny 36}$                  -(-6+aa)/(-1+aa)-(-2+2*aa)/(-1+aa)-((-((-4-aa)^3/(-1+aa)^3)...$\\
 \hspace*{-2em}{\tiny 37}$                  +(4*(-4-aa)*(-6+aa))/(-1+aa)^2-(8*(-4-aa))/(-1+aa))*sqrt(1-2*aa+aa^2))...$\\
 \hspace*{-2em}{\tiny 38}$                  /(2*sqrt(5)*sqrt(4*aa+aa^2)));$\\
 \hspace*{-2em}{\tiny 39}$                   $\color{mblue}$while$\color{black} $ it<maxiter  abs(c1)>1e-2 && it<maxiter && (abs(c1))<1000$\\
 \hspace*{-2em}{\tiny 40}$                      c1=-c1.^4*(-aa+(1+c1)^4)/(-1-4*c1-6*c1^2-4*c1^3+c1^4*(aa-1));$\\
 \hspace*{-2em}{\tiny 41}$                      it=it+1;$\\
 \hspace*{-2em}{\tiny 42}$                  $\color{mblue}$end$\color{black} \\
 \hspace*{-2em}{\tiny 43}$                  c(j,k)=c1;$\\
 \hspace*{-2em}{\tiny 44}$                  $\color{mblue}$if $\color{black} $abs(c1)<1e-2$\\
 \hspace*{-2em}{\tiny 45}$                         R(j,k)=it/maxiter; $\\
 \hspace*{-2em}{\tiny 46}$                  $\color{mblue}$else if $\color{black} $ abs(c1)>=1000 \|\|  ~isfinite(c1)$\\
 \hspace*{-2em}{\tiny 47}$                      R(j,k)=it/maxiter;$\\
 \hspace*{-2em}{\tiny 48}$                      $\color{mblue}else$\color{black} $ \\
 \hspace*{-2em}{\tiny 49}$                      R(j,k)=1;$\\
 \hspace*{-2em}{\tiny 50}$                      G(j,k)=1;$\\
 \hspace*{-2em}{\tiny 51}$                      B(j,k)=1;$\\
 \hspace*{-2em}{\tiny 52}$                      $\color{mblue}$end$\color{black}  \\
 \hspace*{-2em}{\tiny 53}$                  $\color{mblue}$end$\color{black} \\
 \hspace*{-2em}{\tiny 54}$              $\color{mblue}$end$\color{black} \\
 \hspace*{-2em}{\tiny 55}$          $\color{mblue}$end$\color{black} \\
 \hspace*{-2em}{\tiny 56}$  $\color{mgreen}$\% Image display$\color{black}\\
 \hspace*{-2em}{\tiny 57}$  I(:,:,1)=R(:,:); I(:,:,2)=G(:,:); I(:,:,3)=B(:,:);$\color{black}\\
 \hspace*{-2em}{\tiny 58}$  figure, imshow(I,$\color{mred}$'Xdata'$\color{black}$,[x0 xN], $\color{mred}$'Ydata'$\color{black}$, [y0 yN])$\\
 \hspace*{-2em}{\tiny 59}$  axis $\color{mred}$on$\color{black}$, axis $\color{mred}$xy$\color{black}$, hold $\color{mred}$on$\color{black}\\
 \hspace*{-2em}{\tiny 60}$  plot(real(pa),imag(pa),$\color{mred}$'w*'$\color{black}$)$\\
 \hspace*{-2em}{\tiny 61}$  xlabel($\color{mred}$'Re\{z\}'$\color{black}$); ylabel($\color{mred}$'Im\{z\}'$\color{black}$);$\\
 \hspace*{-2em}{\tiny 62}$  axis xy$

\UndefineShortVerb{\$}
\UndefineShortVerb{\#}

The code is divided into five different parts:

\begin{enumerate}
  \item Values (lines 17--18).\\ The bounds are renamed and the symbolic function introduced as {\ttfamily fun} is translated to an anonymous function, recallable by the output handle.
  \item Fixed point operators (lines 23).\\
  \item Calculation of attractive fixed points (lines 26--36).
  \item Image creation (lines 39--94).\\ Once the fixed point operator and the attracting points are set, the next step consists of the determination of the basins of attraction. The combination of the input parameters {\ttfamily bounds} and {\ttfamily points} set the resolution of the image, and it establishes the mesh of complex points (lines 39--50).

      Lines 58--87 are devoted to assign a color to each starting point. It depends on the basin of attraction and the number of iterations needed to reach the root. If the orbit tends to the attracting point set in the first index of line 35, the point is pictured in orange, as lines 67--69 show; for cases second and third, the point is pictured in blue (lines 72--74) and green (lines 78--80), respectively. Otherwise, the point is not modified, so its color is black.

      As the number of iterations needed to reach convergence increases, its corresponding color gets closer to white (black in the decreasing case). A coefficient in each case (lines 66, 71 and 77) is high if the number of iterations is low, and the RGB values are greater than in the slow orbit instance.
  \item Image display (lines 90--94).\\ The image display is based on the {\ttfamily imshow} command. Images are usually displayed in matrix form (from top to bottom and from left to right). In this case, the image is composed by complex points, so the natural display is the cartesian one (from bottom to top and from left to right). With this purpose, {\ttfamily axis xy} is written in line 95.
\end{enumerate}

Once the program is executed, the output values are the image {\ttfamily I} and the number of iterations of each point {\ttfamily it}. Our recommendation is the use of the {\ttfamily surf} command to plot the number of iterations, in combination with the {\ttfamily shading} one.

In order to apply the introduced code to different fixed point operators, the only part to be changed is the fixed point operators corresponding one. If the method can converge to more than three points, just add another {\ttfamily else if} structure (as lines 79--84) and set a color as many times as necessary.




\begin{thebibliography}{99}
\bibitem{DH} A. Douady, J.H. Hubbard, On the dynamics of polynomials-like mappings, Ann. Sci. Ec. Norm. Sup., 18 (1985) 287-343.
\bibitem{CGS} J. Curry, L. Garnet, D. Sullivan, On the iteration of a rational function: Computer experiments with Newton's method, Comm. Math. Phys., 91 (1983) 267-277.
\bibitem{V} J.L. Varona,
            Graphic and numerical comparison between iterative methods,
            Math. Intelligencer, 24(1) (2002) 37-46.
\bibitem{Amat3} S. Amat, S. Busquier, S. Plaza, Review of some iterative root-finding methods from a dynamical point of view, Scientia Series A: Mathematical Sciences, 10 (2004) 3-35.
\bibitem{Blanchard2}P. Blanchard, The dynamics of Newton's method, Proc. of Symposia in Applied Math., 49 (1994) 139-154.
\bibitem{Fallega} N. Fagella, Invariants in din\`{a}mica complexa, Butllet\'{\i} de la Soc. Cat. de Matem\`{a}tiques, 23(1) (2008) 29-51.
\bibitem{GHR} J.M. Guti\'errez, M.A. Hern\'andez and N. Romero,
              Dynamics of a new family of iterative processes for quadratic polynomials,
              J. of Computational and Applied Mathematics,  233 (2010) 2688-2695.
\bibitem{HPR} G. Honorato, S. Plaza, N. Romero,
            Dynamics of a high-order family of iterative methods,
            Journal of Complexity, 27 (2011) 221-229.

\bibitem{CCGT} F. Chicharro, A. Cordero, J.M. Guti\'errez, J.R. Torregrosa,
            Complex dynamics of derivative-free methods for
            nonlinear equations,
            Applied Mathematics and Computation, doi: 10.1016/j.amc.2012.12.075.
\bibitem{ACCT} S. Artidiello, F. Chicharro, A. Cordero, J.R. Torregrosa,
            Local convergence and dynamical analysis of a new family of optimal fourth-order iterative methods,
            International Journal of Computer Mathematics doi:10.1080/00207160.2012.748900.

\bibitem{SNC} M. Scott, B. Neta, C. Chun,
            Basin attractors for various methods,
            Applied Mathematics and Computation, 218 (2011)
            2584-2599.
\bibitem{CLND} C. Chun, M.Y. Lee, B. Neta, J. D\v{z}uni\'{c},
            On optimal fourth-order iterative methods free from
            second derivative and their dynamics,
            Applied Mathematics and Computation, 218 (2012)
            6427-6438.
\bibitem{NSC1} B. Neta, M. Scott, C. Chun,
            Basin attractors for various methods for multiple roots,
         Applied Mathematics and Computation, 218 (2012) 5043-5066.


%

\bibitem{devaney} R.L. Devaney,
            The Mandelbrot Set, the Farey Tree and the Fibonacci sequence,
            Am. Math. Monthly, 106(4) (1999) 289-302.

\bibitem{kim} Y. I. Kim,
            A triparametric family of three-step optimal
            eighth-order methods for solving nonlinear equations,
            International Journal of Computer Mathematics, 89(8) (2012) 1051-1059.
\bibitem{King_AML} Alicia Cordero, Javier Garcia-Maimo, Juan R. Torregrosa, Maria P. Vassileva, Pura Vindel,
            Chaos in King's iterative family,
            Applied Mathematics Letters, doi: 10.1016/j.aml.2013.03.012.












\end{thebibliography}
\end{document}